\documentclass[12pt,oneside,reqno]{amsart}
\usepackage{amsmath}
\usepackage{graphicx}
\usepackage{mathrsfs}
\usepackage{stmaryrd}
\usepackage{amsfonts}
\usepackage{enumerate,amsmath,amssymb,amsthm}
\numberwithin{equation}{section}

\newcommand{\be}{\begin{eqnarray}}
\newcommand{\ee}{\end{eqnarray}}
\newcommand{\ce}{\begin{eqnarray*}}
\newcommand{\de}{\end{eqnarray*}}
\newtheorem{theorem}{Theorem}[section]
\newtheorem{lemma}[theorem]{Lemma}
\newtheorem{remark}[theorem]{Remark}
\newtheorem{definition}[theorem]{Definition}
\newtheorem{proposition}[theorem]{Proposition}

\newtheorem{corollary}[theorem]{Corollary}

\def\a{\alpha}
\def\eps{\epsilon}

\def\om{\omega}

\def\p{\partial}

\def\l{\lambda}

\def\[{{\Big[}}
\def\]{{\Big]}}
\def\<{{\langle}}
\def\>{{\rangle}}
\def\({{\Big(}}
\def\){{\Big)}}

\def\dif{{\mathord{{\rm d}}}}

\def\div{{\mathord{{\rm div}}}}

\def\u{\mathord{{\bf u}}}

\def\v{\mathord{{\bf v}}}

\def\h{\mathord{{\bf h}}}

\def\no{\nonumber}
\def\bt{\begin{theorem}}
\def\et{\end{theorem}}
\def\bl{\begin{lemma}}
\def\el{\end{lemma}}
\def\br{\begin{remark}}
\def\er{\end{remark}}
\def\bd{\begin{definition}}
\def\ed{\end{definition}}
\def\bp{\begin{proposition}}
\def\ep{\end{proposition}}
\def\bc{\begin{corollary}}
\def\ec{\end{corollary}}

\def\cB{{\mathcal B}}

\def\cF{{\mathcal F}}

\def\cL{{\mathcal L}}

\def\cN{{\mathcal N}}
\def\cO{{\mathcal O}}
\def\cP{{\mathcal P}}

\def\cT{{\mathcal T}}

\def\mE{{\mathbb E}}

\def\mH{{\mathbb H}}

\def\mN{{\mathbb N}}

\def\mP{{\mathbb P}}

\def\mR{{\mathbb R}}

\def\mW{{\mathbb W}}

\def\mZ{{\mathbb Z}}

\def\geq{\geqslant}
\def\leq{\leqslant}

\def\sB{{\mathscr B}}

\def\sD{{\mathscr D}}

\def\sP{{\mathscr P}}

\allowdisplaybreaks

\begin{document}

\bf\title{Smooth Solutions of Non-linear Stochastic
Partial Differential Equations}
\author{Xicheng Zhang}
\dedicatory{Department of Mathematics,
Huazhong University of Science and Technology\\
Wuhan, Hubei 430074, P.R.China\\
Email: XichengZhang@gmail.com}
\begin{abstract}

In this paper, we study the regularities of solutions of
nonlinear stochastic partial differential equations in the framework of Hilbert scales. Then we apply our general result to
several typical nonlinear SPDEs such as stochastic Burgers and Ginzburg-Landau's  equations
on the real line, stochastic 2D  Navier-Stokes equations in the whole space
and a stochastic tamed 3D Navier-Stokes equation in the whole space, and
obtain the existence of their respectively smooth solutions.
\end{abstract}

\maketitle

\rm
\tableofcontents
\section{Introduction}

Consider the following stochastic Burgers and Ginzburg-Landau equation on the real line:
\be
\label{BGL}\left\{
\begin{array}{lcl}
\dif u(t,x)&=&\Big[\nu\p^2_x u(t,x)+c_0\cdot\p_x u(t,x)^2+c_1\cdot u(t,x)-c_2\cdot u(t,x)^3\Big]\dif t\\
&&+\sum_k\sigma_k(t,x,u(t,x))\dif W^k(t),\\
u(0,x)&=&u_0(x),\ \ \ x\in\mR,
\end{array}
\right.
\ee
where $c_0,c_1\in\mR$ and $\nu,c_2>0$, $\{W^k(t),k\in\mN\}$ is a sequence
of independent Brownian motions,
the coefficients $\{\sigma_k,k\in\mN\}$ satisfy some smoothness conditions.
Up to now,  there are many papers devoted to the studies of stochastic Burgers'
equation and
stochastic Ginzburg-Landau's equation (cf. \cite{Be-Ca, Gy-Nu, Ro-So,Ki} and references therein).
In \cite{Gy-Nu}, using heat kernel estimates,
Gy\"ongy and Nualart proved the existence and uniqueness
of $L^2(\mR)$-solution to stochastic Burgers' equation on the real line.
By solving an infinite dimensional Kolmogorov's equation, R\"ockner and Sobol \cite{Ro-So}
developed a new method to solve the generalized stochastic Burgers and reaction diffusion equations.
More recently, Kim \cite{Ki} studied the stochastic Burgers type equation
with the first order term having polynomial growth, as well as the
existence of associated invariant measures.

Since all of these works are concerned with stochastic Burgers equation
driven by space-time white noises, they had to consider weak or mild solutions
rather than strong or classical solutions.
A natural question is that: does there exist smooth solution or classical solution  in $x$
to the equation (\ref{BGL}) if all the datas are smooth? Of course,
for this question, we can only consider the
equation (\ref{BGL}) driven by time white and space correlated noises.
We remark that for the deterministic Burgers equation, i.e., $\sigma_k=c_1=c_2=0$
and $c_0=1, \nu>0$, it is well known that there
exists a unique smooth solution if the initial data is smooth  (cf. \cite{Kre}).

Let us also consider the following stochastic 2D Navier-Stokes equation
in $\mR^2$:
\be
\label{Ns3}\left\{
\begin{array}{lcl}
\p_t u_1=\nu\Delta u_1-u_1\p_{x_1}u_1-u_2\p_{x_2}u_1-\p_{x_1}p+f_1\\
\p_t u_2=\nu\Delta u_2-u_1\p_{x_1}u_2-u_2\p_{x_2}u_2-\p_{x_2}p+f_2\\
\p_{x_1} u_1+\p_{x_2}u_2=0,
\end{array}
\right.
\ee
where $\nu$ is the viscosity constant, $\u(t,x)=(u_1,u_2)$ is the velocity field,
$p$ is the pressure function, and $f=(f_1,f_2)$ is the white in time
and additive stochastic forcing. In \cite{MiRo}, Mikulevicius and Rozovskii studied
the existence of martingale solutions for any dimensional stochastic Navier-Stokes equations
in the whole space. In particular, they obtained the
existence of a unique weak solution for the above two dimensional equation.
In the periodic boundary case,
using Galerkin's approximation
and Fourier's transformation, Mattingly \cite{Ma} proved the spatial analyticity
for the solution to the above stochastically forced 2D Navier-Stokes equation.
However, using his method, it seems to be hard to consider
the multiplicative noise force.

As for the stochastic 3D Navier-Stokes equation, R\"ockner and the author
\cite{Ro-Zh, Ro-Zh1} recently studied the following tamed or modified scheme
in the whole space $\mR^3$ or periodic boundary case:
\be
\label{Ns4}\left\{
\begin{array}{lcl}
\p_t u_j=\nu\Delta u_j-\sum_{i=1}^3u_i\p_{x_i}u_j-\p_{x_j}p\\
\qquad\quad -g_N\Big(\sum_{i=1}^3|u_i|^2\Big)u_j+f_j,\ \ \ \ j=1,2,3,\\
\sum_{i=1}^3\p_{x_i} u_i=0,
\end{array}
\right.
\ee
where the taming function $g_N:\mR_+\mapsto\mR_+$ is smooth and satisfies that
\ce
g_N(r)=0 \mbox{ on } r\leq N\mbox{ and } g_N(r)=(r-N)/\nu \mbox{ on } r\geq N+2.
\de
In \cite{Ro-Zh1}, we proved the existence of a unique strong solution and
the ergodicity of associated invariant measure in the case of periodic boundary conditions.
For the existence, the method is mainly based on the Galerkin approximation, and the smooth solution of
Eq.(\ref{Ns4}) is not obtained.

Our main purpose in this paper is to present a unified settings for
proving the existence of smooth solutions to the above three typical nonlinear stochastic
partial differential equations. That is, we shall consider an abstract semilinear
stochastic evolution equation  in the scope of Hilbert scales
determined by a sectorial operator.
Here, the analytic semigroup generated by the sectorial operator plays a mollifying
role, and will be used to construct an approximating sequence of
regularized stochastic ordinary differential equations in Hilbert spaces.
After obtaining some uniform estimates of the approximating solutions
in the spaces of Hilbert scale,
we can prove that the solutions of approximating equations strongly converge
to a smooth solution. Our approach is much inspired by the energy
method used in the deterministic case (cf. \cite{Ma-Ber}), and is different
from Galerkin's approximation and semigroup methods which were extensively
used in the well known studies of stochastic partial differential equations(abbrev. SPDEs)
 (cf. \cite{Kr-Ro, DaZa}, etc.).
We remark that  the regularity of solution will be decreasing
when we use the semigroup method to deal with  SPDEs (cf. \cite[Sections
5,8]{b1}). The main advantage of our method is that we can obtain
better regularities unlike the semigroup method.

In \cite{Zh}, using the semigroup method and nonlinear interpolations,
we have already proved the existence of smooth solutions to a large class
of semilinear SPDEs when the coefficients are smooth and have all bounded derivatives.
However, the result in \cite{Zh} can not be applied to the above mentioned equations.
It should be emphasized that the existence of smooth solutions
for nonlinear partial differential equations, fox example, Navier-Stokes equations,
usually depends on the spatial dimension. Thus, it is not expected to use
our general result(see Theorem \ref{main} below)
to treat  high dimensional nonlinear SPDEs for obtaining smooth solutions.
Nevertheless, we may still apply our general result to achieving the existence of strong
solutions for a class of   semilinear SPDEs with Lipschitz nonlinear coefficients
in Euclidean space (cf. \cite{Kr1, Kr2, MiRo}).
We also want to say that although our main attention concentrates
on the above three types nonlinear SPDEs,
our result can also be applied to dealing with the stochastic Kuramoto-Sivashinsky equation
and stochastic Cahn-Hilliard equation (cf. \cite{Se-You}),
as well as the stochastic partial differential equation in the abstract
Wiener space (cf. \cite{Zh2}), which are not discussed here.

This paper is organized as follows:
in Section 2, we shall give the general framework and state our main result.
In Section 3, we devote to the proof of our main result. In Section 4,
we investigate a class of semilinear SPDEs in the whole space and in bounded smooth domains
of Euclidean space,
and obtain the existence of unique strong solutions. In Section 5,
we study stochastic Burgers and Ginzburg-Landau's equations on the real line, and get
the existence of smooth solutions.
In Section 6, we prove the existence of smooth solutions to stochastic tamed 3D Navier-Stokes
equations. In Section 7, stochastic 2D Navier-Stokes equation  in the whole space
and with multiplicative noise
is considered.

\textsc{Convention}: The letter $C$ with or without subscripts will
denote a positive constant, which is unimportant and may change from
one line to another line.

\section{General Settings and Main Result}

Let $(\mH,\|\cdot\|_\mH)$ be a separable Hilbert space, $\cL$ a symmetric and non-positive
sectorial operator in $\mH$ that generates a symmetric analytic semigroup
$(\cT_\eps)_{\eps>0}$ in $\mH$ (cf. \cite{Pa}).
For $\a\geq 0$, we define the Sobolev space $\mH^\a$ by
$$
\mH^\a:=\sD((I-\cL)^{\a/2})
$$
together with the norm
$$
\|u\|_\a:=\|(I-\cL)^{\a/2}u\|_{\mH}.
$$
The inner product in $\mH^\a$ is denoted by $\<\cdot,\cdot\>_\a$.
The dual space of $\mH^\a$ is denoted by $\mH^{-\a}$ with the norm
$$
\|u\|_{-\a}:=\|(I-\cL)^{-\a/2}u\|_{\mH}.
$$
Then $(\mH^\a)_{\a\in\mR}$ forms a Hilbert scale (cf. \cite{Kr-Re-Se,Ro}), i.e.:
\begin{enumerate}[(i)]
\item for any $\a<\beta$, $\mH^\beta\subset\mH^\a$;
\item for any $\a<\gamma<\beta$ and $u\in\mH^\beta$,
\be
\|u\|_\gamma\leq C\|u\|^{\frac{\beta-\gamma}{\beta-\a}}_\a
\cdot\|u\|_\beta^{\frac{\gamma-\a}{\beta-\a}}.\label{Int}
\ee
\end{enumerate}

Set $\mH^\infty:=\cap_{m\in\mN}\mH^m$. Then  (cf. \cite{Pa})
\bp\label{Pro1}
For all integer $m\in\mZ$, we have
\begin{enumerate}[(i)]
\item $\mH^\infty$ is dense in $\mH^m$, and for every $\eps>0$ and $u\in\mH^m$,
$\cT_\eps u\in\mH^\infty$;
\item for every  $\eps>0$ and $u\in\mH^m$
$$
(I-\cL)^{m/2}\cT_\eps u=\cT_\eps (I-\cL)^{m/2}u;
$$
\item for every  $\eps>0$ and $u\in\mH^m$, $k=1,2,\cdots$
$$
\|\cT_\eps u-u\|_m\leq C_m\eps^{k/2}\|u\|_{m+k};
$$
\item for every  $\eps>0$ and $u\in\mH^m$, $k=0,1,2,\cdots$
$$
\|\cT_\eps u\|_{m+k}\leq \frac{C_{mk}}{\eps^{k/2}}\|u\|_m.
$$
\end{enumerate}
\ep
Let $l^2$ be the usual Hilbert space of square summable sequences of real numbers.
Let $(\Omega,\cF,(\cF_t)_{t\geq 0}, \mP)$ be a complete filtration probability space.
A family of independent one dimensional $\cF_t$-adapted Brownian motions
$\{W^k(t); t\geq 0, k=1,2,\cdots\}$ on $(\Omega,\cF, \mP)$ are given. Then $\{W(t), t\geq 0\}$
can be regarded as a cylindrical Brownian motion in $l^2$ (cf. \cite{DaZa}).

Consider the following type stochastic evolution equation
\be
\dif u(t)=[\cL u(t)+F(t,u(t))]\dif t+\sum_kB_k(t,u(t))\dif W^k(t),
\quad u(0)=u_0,
\ee
where the stochastic integral is understood as It\^o's integral,
and for some $N\in\mN$ the coefficients
\ce
F(t,\om,u): &&\mR_+\times\Omega\times\mH^N\to\mH^{-1},\\
B(t,\om,u): &&\mR_+\times\Omega\times\mH^0\to\mH^0\otimes l^2
\de
are two measurable functions, and for every $t\geq0$ and $u\in\mH^N$,
$$
F(t,\cdot,u)\in\cF_t/\sB(\mH^{-1}),\ \ B(t,\cdot,u)\in\cF_t/\sB(\mH^0\otimes l^2).
$$
We also require that $F(t,\om,u)\in\mH^0$ for $u\in\mH^{N+1}$
and $B_k(t,\om,u)\in\mH^m$ for any $m\in\mN$ and $u\in\mH^{m+1}$.

We make the following assumptions on $F$ and $B$:
\begin{enumerate}[{(\bf $\mathbf{H}1_N$})]
\item
There exist $q_1,q_2\geq 1$ and constants $\l_0,\l_1,\l_2,\l_3,\l_4>0$ such that
for all $(t,\om)\in\mR_+\times\Omega$
and $u,v\in\mH^N$
\be
\|F(t,\om,u)-F(t,\om,v)\|_{-1}&\leq& \l_0(\|u\|^{q_1}_N+\|v\|^{q_1}_N+1)\cdot\|u-v\|_0,\label{Co1}\\
\sum_k\|B_k(t,\om,u)-B_k(t,\om,v)\|_0^2&\leq& \l_1\|u-v\|^2_0,\label{Co3}
\ee
and for $u\in\mH^{N+1}$
\be
\<u,F(t,\om,u)\>_0&\leq& \frac{1}{2}\|u\|_1^2+\l_2(\|u\|_0^2+1),\label{Co2}\\
\|F(t,\om,u)\|_0&\leq&\l_3(\|u\|_{N+1}+\|u\|^{q_2}_N+1),\label{Co7}\\
\sum_k\|B_k(t,\om,u)\|^2_0&\leq&\l_4(\|u\|^2_0+1).\label{Co9}
\ee
\end{enumerate}
\begin{enumerate}[{(\bf $\mathbf{H}2_\cN$})]
\item
For some integer $\cN\geq N$,  and each $m=1,\cdots, \cN$, any $\delta\in(0,1)$,
there exist $\a_m,\beta_m\geq 1$
and constants $\l_{1m},\l_{2m}>0$
such that for all $u\in\mH^\infty$
and $(t,\om)\in\mR_+\times\Omega$
\be
\<u,F(t,\om,u)\>_m&=&\<(I-\cL)^{m+\frac{1}{2}}u,(I-\cL)^{-\frac{1}{2}}F(t,\om,u)\>_0\no\\
&\leq& \frac{1}{2}\|u\|_{m+1}^2+\l_{1m}(\|u\|^{\a_m}_{m-1}+1),\label{Co4}\\
\sum_k\|B_k(t,\om,u)\|_m^2&\leq& \delta\|u\|^2_{m+1}+\l_{2m}(\|u\|^{\beta_m}_{m-1}+1).\label{Co6}
\ee
\end{enumerate}

Our main result is that
\bt\label{main}
Under {(\bf $\mathbf{H}1_N$}) and {(\bf $\mathbf{H}2_\cN$}) with $\cN\geq N\geq 1$,
for any $u_0\in\mH^\cN$, there exists a unique process
$u(t)\in\mH^\cN$ such that
\begin{enumerate}[(i)]
\item The process $t\mapsto u(t)\in\mH^0$ is $\cF_t$-adapted and continuous, and
for any $T>0$ and $p\geq 2$
$$
\mE\left(\sup_{s\in[0,T]}\|u(s)\|^p_\cN\right)
+\int^T_0\mE\|u(s)\|^2_{\cN+1}\dif s<+\infty.
$$
\item $u(t)$ satisfies the following equation in $\mH^0$: for all $t\geq 0$
$$
u(t)=u_0+\int^t_0[\cL u(s)+F(s,u(s))]\dif s+\sum_k\int^t_0B_k(s,u(s))\dif W^k(s), \ \ \mP-a.s..
$$
\end{enumerate}
\et
\br
By (\ref{Co7}), (\ref{Co9}) and (i), one knows that all the integrals appearing
in (ii) are well defined. Moreover, if for some $C>0$, $p\geq 1$ and any $u\in\mH^{\cN+1}$
$$
\|F(s,u)\|_{\cN-1}\leq C(\|u\|_{\cN+1}+\|u\|^{p}_\cN+1),
$$
then we can find an $\mH^\cN$-valued continuous version of $u$  (cf. \cite{Ro}).
\er
\br
The solution $u(t)$ also satisfies the following integral equation written in terms of
the semigroup $\cT_t$:
$$
u(t)=\cT_t u_0+\int^t_0\cT_{t-s}F(s,u(s))\dif s+\sum_k\int^t_0\cT_{t-s}B_k(s,u(s))\dif W^k(s).
$$
Using this representation, we can further prove the H\"older continuity of $u(t)$
in $t$  (cf. \cite{Zh2}).
\er
\section{Proof of Main Theorem}

\subsection{Regularized Stochastic Differential Equations}
For any $m=0,\cdots,\cN$, let us consider the following regularized
stochastic differential equation in $\mH^m$:
\be
\dif u^\eps(t)=A^\eps(t,u^\eps(t))\dif t+\sum_kB_k^\eps(t,u^\eps(t))\dif W^k(t),\ \  u^\eps(0)=u_0,
\label{Es5}
\ee
where the regularized operators are given by:
\ce
A^\eps(t,\om,u)&:=&\cT_\eps\cL\cT_\eps u+\cT_\eps F(t,\om,\cT_\eps u)\\
B^\eps_k(t,\om,u)&:=&\cT_{\eps} B_k(t,\om,\cT_\eps u).
\de
The following two lemmas are direct from {(\bf $\mathbf{H}1_N$}) and {(\bf $\mathbf{H}2_\cN$}).
We omit the proof.
\bl\label{Io4}
There exists a constant $C>0$ such that for any $\eps>0$ and all
$(t,\om)\in\mR_+\times\Omega$, $u\in\mH^0$
\ce
\<u,A^\eps(t,\om,u)\>_0&\leq& -\frac{1}{2}\|\cT_\eps u\|^2_1+C(\|u\|^2_0+1),\\
\sum_k\|B^\eps_k(t,\om,u)\|_0^2&\leq&C(\|u\|^2_0+1).
\de
\el
\bl\label{Le6}
For any $m=1,\cdots,\cN$ and $\delta\in(0,1)$,
there exist two constants $C_m,C_{m,\delta}>0$ such that
for any $\eps>0$ and all $(t,\om)\in\mR_+\times\Omega$ and $u\in\mH^m$, we have
\ce
\<u,A^\eps(t,\om,u)\>_m&\leq& -\frac{1}{2}\|\cT_\eps u\|^2_{m+1}
+C_m(\|u\|^{\a_m}_{m-1}+1)\\
\sum_k\|B^\eps_k(t,\om,u)\|_m^2&\leq&\delta\|\cT_\eps u\|^2_{m+1}
+C_{m,\delta}(\|u\|^{\beta_m}_{m-1}+1),
\de
where $\a_m$ and $\beta_m$ are same as in {(\bf $\mathbf{H}2_\cN$}).
\el
We also have
\bl\label{Es1}
For any $m=0,\cdots,\cN$, the functions $A^\eps$ and $B^\eps$
are locally Lipschitz continuous in $\mH^m$ with respect to $ u$.
More precisely, for any $R>0$ there are $C_{R,\eps}, C'_{R,\eps}>0$
such that for any $(t,\om)\in\mR_+\times\Omega$ and $u, v\in \mH^m$
with $\| u\|_m,\| v\|_m\leq R$
\ce
\|A^\eps(t,\om, u)-A^\eps(t,\om, v)\|_m&\leq& C_{R,\eps}\| u- v\|_m\\
\sum_k\|B^\eps_k(t,\om, u)-B^\eps_k(t,\om, v)\|^2_m&\leq& C'_{R,\eps}\| u- v\|_m^2.
\de
\el
\begin{proof}
By (iv) of Proposition \ref{Pro1}, we have
$$
\|\cT_\eps \cL(\cT_\eps u)-\cT_\eps \cL(\cT_\eps v)\|_m
=\|\cL\cT_\eps^2( u- v)\|_{m}\leq C_\eps\| u- v\|_{m},
$$
and by {(\bf $\mathbf{H}1_N$})
\ce
&&\|\cT_\eps F(t,\cT_\eps u)-\cT_\eps F(t,\cT_\eps v)\|^2_m
+\sum_k\|\cT_\eps B_k(t,\cT_\eps u)-\cT_\eps B_k(t,\cT_\eps v)\|_m^2\\
&\leq&C_\eps\|F(t,\cT_\eps u)-F(t,\cT_\eps v)\|_{-1}^2
+C_\eps\sum_k\|B_k(t,\cT_\eps u)-B_k(t,\cT_\eps v)\|_0^2\\
&\leq& C_{R,\eps}\|\cT_\eps u- \cT_\eps v\|_{0}^2\leq C_{R,\eps}\|u- v\|_m^2.
\de
The proof is complete.
\end{proof}

We now prove the following key estimate about the solution of regularized
stochastic differential equation (\ref{Es5}).
\bt\label{Th1}
For any $u_0\in\mH^\cN$, there exists a unique continuous $\cF_t$-adapted
solution $u^\eps$ to Eq.(\ref{Es5}) in $\mH^\cN$
such that for any $p\geq 1$ and $T>0$
\be
\sup_{\eps\in(0,1)}\mE\left(\sup_{t\in[0,T]}\|u^\eps(t)\|^{2p}_\cN\right)
+\sup_{\eps\in(0,1)}\int^T_0\mE\|\cT_\eps u^\eps(s)\|^2_{\cN+1}\dif s\leq C_{p,T}.\label{Es2}
\ee
\et
\begin{proof}
First of all, by Lemma \ref{Es1},
there exists a unique continuous $\cF_t$-adapted
local solution $u^\eps(t)$ in $\mH^m$ for any $m=0,\cdots, \cN$.
We now use induction method to prove
that
\ce
(\sP_m): \left\{
\begin{array}{ll}
\mbox{$u^\eps(t)$ is non-explosive in $\mH^m$ and for any $p\geq 1$ and $T>0$}\\
\sup_{\eps\in(0,1)}\mE\left(\sup_{t\in[0,T]}\|u^\eps(t)\|^{2p}_m\right)\leq C_{m,p,T}.
\end{array}
\right.
\de

By the standard stopping times technique, it suffices to prove the estimate in
$(\sP_m)$. In the following, we fix $T>0$. By It\^o's formula, we have for any $p\geq 1$
\be
\|u^\eps(t)\|^{2p}_m=\|u_0\|^{2p}_m+I_{m1}(t)+I_{m2}(t)+I_{m3}(t)+I_{m4}(t),\label{Io1}
\ee
where
\ce
I_{m1}(t)&:=&2p\int^t_0\|u^\eps(s)\|^{2(p-1)}_m\<u^\eps(s),A^\eps(u^\eps(s))\>_m\dif s\\
I_{m2}(t)&:=&2p\sum_{k=1}^\infty\int^t_0\|u^\eps(s)\|^{2(p-1)}_m\<u^\eps(s),
B^\eps_k(s,u^\eps(s))\>_m\dif W^k_s\\
I_{m3}(t)&:=&p\sum_{k=1}^\infty\int^t_0\|u^\eps(s)\|^{2(p-1)}_m
\|B^\eps_k(s,u^\eps(s))\|_m^2\dif s\\
I_{m4}(t)&:=&2p(p-1)\sum_{k=1}^\infty\int^t_0\|u^\eps(s)\|^{2(p-2)}_m|\<u^\eps(s),
B^\eps_k(s,u^\eps(s))\>_m|^2\dif s.
\de

Let us first look at the case of $m=0$. By Lemma \ref{Io4} and Young's inequality, we have
\be
I_{01}(t)+I_{03}(t)+I_{04}(t)\leq C_p\int^t_0(\|u^\eps(s)\|^{2p}_0+1)\dif s.\label{Io2}
\ee
Taking expectations for (\ref{Io1}) gives that
$$
\mE\|u^\eps(t)\|^{2p}_0\leq \|u_0\|^{2p}_0+C_p\int^t_0(\mE\|u^\eps(s)\|^{2p}_0+1)\dif s.
$$
By Gronwall's inequality, we obtain for any $p\geq 1$
\be
\sup_{t\in[0,T]}\mE\|u^\eps(t)\|^{2p}_0 \leq C_{p,T}(\|u_0\|^{2p}_0+1).\label{Io3}
\ee

On the other hand, by BDG's inequality and Lemma \ref{Io4}, we have
\ce
\mE\left(\sup_{s\in[0,T]}|I_{02}(s)|\right)&\leq&C_p
\mE\left(\int^T_0\|u^\eps(s)\|^{4(p-1)}_0\|\<u^\eps(s),
B^\eps_\cdot(s,u^\eps(s))\>_0\|^2_{l^2}\dif s\right)^{1/2}\\
&\leq&C_p\mE\left(\int^T_0\|u^\eps(s)\|^{4p-2}_0
\cdot\Big(\sum_k\|B^\eps_k(s,u^\eps(s))\|^2_0\Big)\dif s\right)^{1/2}\\
&\leq&C_{p}\mE\left(\int^T_0(\|u^\eps(s)\|^{4p}_0+1)\dif s\right)^{1/2}\\
&\leq&C_{p}\left(\int^T_0(\mE\|u^\eps(s)\|^{4p}_0+1)\dif s\right)^{1/2}\\
&\leq& C_{p,T}(\|u_0\|^{2p}_0+1).
\de
Thus, from (\ref{Io1})-(\ref{Io3}), one knows that $(\sP_0)$ holds.

Suppose now that $(\sP_{m-1})$ holds.
By Lemma \ref{Le6} and Young's inequality, we have
\be
I_{m1}(t)&\leq& p\int^t_0\Big[-\|u^\eps(s)\|^{2(p-1)}_{m}
\|\cT_\eps u^\eps(s)\|^{2}_{m+1}\no\\
&&+C_m\|u^\eps(s)\|^{2(p-1)}_m\cdot(\|u^\eps(s)\|^{\a_m}_{m-1}+1)\Big]\dif s\no\\
&\leq&-p\int^t_0\|u^\eps(s)\|^{2(p-1)}_{m}\|\cT_\eps u^\eps(s)\|^{2}_{m+1}\dif s+
C_{m,p}\int^t_0\|u^\eps(s)\|^{2p}_m\dif s\no\\
&&+C_{m,p}\int^t_0(\|u^\eps(s)\|^{p\cdot\a_m}_{m-1}+1)\dif s,\label{Lk1}
\ee
and for any $\delta\in(0,1)$
\be
I_{m3}(t)+I_{m4}(t)&\leq& p(2p-1)\int^t_0\|u^\eps(s)\|^{2(p-1)}_m\cdot
\Big(\sum_k\|B^\eps_k(s,u^\eps(s))\|^2_{m}\Big)\dif s\no\\
&\leq& p(2p-1)\delta\int^t_0\|u^\eps(s)\|^{2(p-1)}_{m}\|\cT_\eps u^\eps(s)\|^{2}_{m+1}\dif s\no\\
&&+C_{m,p}\int^t_0\|u^\eps(s)\|^{2p}_m\dif s
+C_{m,p}\int^t_0(\|u^\eps(s)\|^{p\cdot\beta_m}_{m-1}+1)\dif s.\label{Lk2}
\ee
Choosing $\delta=\frac{1}{2(2p-1)}$ and taking expectations for (\ref{Io1}) gives that
\ce
&&\mE\|u^\eps(t)\|^{2p}_m+\frac{p}{2}\int^t_0\mE\Big(\|u^\eps(s)\|^{2(p-1)}_m
\cdot\|\cT_\eps u^\eps(s)\|^2_{m+1}\Big)\dif s\\
&\leq&\|u_0\|^{2p}_m+C_{m,p}\int^t_0\mE\|u^\eps(s)\|^{2p}_m\dif s
+C_{m,p}\int^t_0\mE(\|u^\eps(s)\|^{p_m}_{m-1}+1)\dif s,
\de
where $p_m:=p\cdot(\a_m\vee\beta_m)$.

By Gronwall's inequality again and $(\sP_{m-1})$, we get for any $p\geq 1$
\be
\sup_{t\in[0,T]}\mE\|u^\eps(t)\|^{2p}_m+
\int^T_0\mE\Big(\|u^\eps(s)\|^{2(p-1)}_{m}\|\cT_\eps u^\eps(s)\|^{2}_{m+1}\Big)
\dif s\leq C_{m,p,T}.\label{Lk3}
\ee
Furthermore, by BDG's inequality and (\ref{Lk3}), we have
\ce
\mE\left(\sup_{s\in[0,T]}|I_{m2}(s)|\right)&\leq&C_p
\mE\left(\int^T_0\|u^\eps(s)\|^{4(p-1)}_m\cdot\|\<u^\eps(s),
B^\eps_\cdot(s,u^\eps(s))\>_m\|^2_{l^2}\dif s\right)^{1/2}\\
&\leq&C_p\mE\left(\int^T_0\|u^\eps(s)\|^{4p-2}_m
\Big(\sum_k\|B^\eps_k(s,u^\eps(s))\|_m^2\Big)\dif s\right)^{1/2}\\
&\leq&C_p\mE\Bigg(\sup_{s\in[0,T]}\|u^\eps(s)\|^{p}_m\cdot\Bigg[
\int^T_0\|u^\eps(s)\|^{2(p-1)}_m\times\\
&&\Big(\delta\|\cT_\eps u^\eps(s)\|^2_{m+1}+C_{m,\delta}(\|u^\eps(s)\|^{\beta_m}_{m-1}+1)
\Big)\dif s\Bigg]^{1/2}\Bigg)\\
&\leq&C_{p}\delta\int^T_0\mE\Big(\|u^\eps(s)\|^{2(p-1)}_m
\|\cT_\eps u^\eps(s)\|^2_{m+1}\Big)\dif s\\
&&+C_{m,p}\int^T_0\mE\Big(\|u^\eps(s)\|^{2(p-1)}_m
(\|u^\eps(s)\|^{\beta_m}_{m-1}+1)\Big)\dif s\\
&&+\frac{1}{2}\mE\left(\sup_{t\in[0,T]}\|u^\eps(t)\|^{2p}_m\right)\\
&\leq&C_{m,p,T}+\frac{1}{2}\mE\left(\sup_{t\in[0,T]}\|u^\eps(t)\|^{2p}_m\right),
\de
which together with (\ref{Io1}) and (\ref{Lk1})-(\ref{Lk3}) yields $(\sP_{m})$.
The proof is complete.
\end{proof}
\subsection{Convergence of $u^\eps(t)$}
\bl\label{Le1}
For any $R>0$, there exists a constant $C_R>0$ such that for any $0<\eps'<\eps<1$,
$(t,\om)\in\mR_+\times\Omega$ and $ u, v\in\mH^{N+1}$ with $\|u\|_{N},\|v\|_{N}\leq R$
\ce
\< u- v, A^\eps(t,\om, u)-A^{\eps'}(t,\om, v)\>_0&\leq&
C_R\cdot\sqrt{\eps}\cdot(1+\|\cT_{\eps'}v\|_{N+1})+C_R\cdot\|u-v\|_0^2,\\
\sum_k\|B^\eps_k(t,\om, u)-B^{\eps'}_k(t,\om, v)\|^2_0&\leq&
C_R\cdot\sqrt{\eps}\cdot(1+\|\cT_{\eps'}v\|^2_2)+C_R\cdot\|u-v\|_0^2.
\de
\el
\begin{proof}
We only prove the first estimate. The second one is similar.

Above of all, by (ii) and (iii) of Proposition \ref{Pro1}, we have
\be
&&\< u- v, \cT_\eps \cL(\cT_\eps u)-\cT_{\eps'} \cL(\cT_{\eps'}  v)\>_0\no\\
&=&\<\cT_\eps(u-v), \cL\cT_{\eps}(u-v)\>_0
+\<u- v, (\cT_{2\eps}-\cT_{2\eps'})\cL v\>_0\no\\
&\leq&-\|\cT_\eps (u-v)\|_1^2+\|\cT_\eps (u-v)\|_0^2
+C_R\cdot\|(\cT_{2\eps}-\cT_{2\eps'})v\|_1\no\\
&\leq&-\|\cT_\eps (u-v)\|_1^2+C\|u-v\|_0^2+
C_R\cdot\sqrt{\eps}\cdot \|\cT_{\eps'}v\|_2.\label{PP2}
\ee
Secondly, we  decompose the term involving $F$ in $A^\eps$ as:
\ce
&&\< u- v, \cT_\eps F(t,\cT_\eps u)-\cT_{\eps'} F(t,\cT_{\eps'} v)\>_0\\
&=&\<\cT_{\eps}(u-v), F(t,\cT_{\eps} u)-F(t,\cT_{\eps'} v)\>_0\\
&&+\<(\cT_{\eps}-\cT_{\eps'})(u-v), F(t,\cT_{\eps'} v)\>_0\\
&=:&I_1+I_2.
\de
By (\ref{Co1}) and (iii) of Proposition \ref{Pro1}, we have  for $I_1$
\ce
I_1&\leq&\frac{1}{4}\|\cT_\eps(u-v)\|_1^2+\|F(t,\cT_{\eps} u)-F(t,\cT_{\eps'} v)\|^2_{-1}\\
&\leq&\frac{1}{4}\|\cT_\eps(u-v)\|_1^2+C_R\cdot\|\cT_{\eps} u-\cT_{\eps'} v\|^2_0\\
&\leq& \frac{1}{4}\|\cT_\eps(u-v)\|_1^2+C_R\cdot\sqrt\eps+
C_R\cdot\|u-v\|_0^2,
\de
and for $I_2$, by (\ref{Co7})
\ce
I_2&\leq& \|(\cT_{\eps}-\cT_{\eps'})(u-v)\|_0\cdot\|F(t,\cT_{\eps'} v)\|_0\\
&\leq&C_R\cdot\sqrt{\eps}\cdot(\|\cT_{\eps'}v\|_{N+1}+1).
\de
Combining the above calculations yields the first one.
\end{proof}

We now prove that
\bl\label{Le7}
For any $T>0$, we have
\ce
\lim_{\eps,\eps'\downarrow 0}\mE\left(\sup_{t\in[0,T]}\|u^\eps(t)-u^{\eps'}(t)\|^{2}_0\right)=0.
\de
\el
\begin{proof}
For any $R>0$ and $1>\eps>\eps'>0$, define the stopping time
\ce
\tau^{\eps,\eps'}_R:=\inf\{t>0: \|u^\eps(t)\|_{N}\wedge\|u^{\eps'}(t)\|_{N}\geq R\}.
\de
Then, by Theorem \ref{Th1} we have
\be
\mP(\tau^{\eps,\eps'}_R<T)\leq\frac{\mE\Big(\sup_{t\in[0,T]}
(\|u^\eps(t)\|^2_{N}\wedge\|u^{\eps'}(t)\|^2_{N})\Big)}{R^2}\leq \frac{C_{T}}{R^{2}}.\label{Io6}
\ee
Set
$$
v(t):=u^\eps(t)-u^{\eps'}(t).
$$
By It\^o's formula, we have
\be
\|v(t)\|^{2}_0=J_1(t)+J_2(t)+J_3(t),\label{Io5}
\ee
where
\ce
J_1(t)&:=&2\int^t_0\<v(s),A^\eps(s,u^\eps(s))-A^{\eps'}(s,u^{\eps'}(s))\>_0\dif s\\
J_2(t)&:=&\sum_k\int^t_0\|B^\eps_k(s,u^\eps(s))-B^{\eps'}_k(s,u^{\eps'}(s))\|^2_0\dif s\\
J_3(t)&:=&2\sum_k\int^t_0\<v(s),B^\eps_k(s,u^\eps(s))
-B^{\eps'}_k(s,u^{\eps'}(s))\>_0\dif W^k_s.
\de
By Lemma \ref{Le1}, we have
\be
&&J_1(t\wedge\tau^{\eps,\eps'}_R)+J_2(t\wedge\tau^{\eps,\eps'}_R)\no\\
&\leq&C_{R}\cdot\sqrt{\eps}\cdot\left(1+\int^T_0\|\cT_{\eps'}u^{\eps'}(s)\|^2_{N+1}\dif s\right)
+C_R\int_0^{t\wedge\tau^{\eps,\eps'}_R}\|v(s)\|^{2}_0\dif s\no\\
&\leq&C_{R}\cdot\sqrt{\eps}\cdot\left(1+\int^T_0\|\cT_{\eps'}u^{\eps'}(s)\|^2_{N+1}\dif s\right)
+C_R\int_0^t\|v(s\wedge\tau^{\eps,\eps'}_R)\|^{2}_0\dif s.\label{PP7}
\ee
Hence,  by Theorem \ref{Th1} and taking expectations for (\ref{Io5}), we obtain
\ce
\mE\|v(t\wedge\tau^{\eps,\eps'}_R)\|^2_0
\leq C_{R,T}\cdot\sqrt{\eps}+C_R\int_0^t\mE\|v(s\wedge\tau^{\eps,\eps'}_R)\|^{2}_0\dif s.
\de
By Gronwall's inequality, we get
\be
\sup_{t\in[0,T]}\mE\|v(t\wedge\tau^{\eps,\eps'}_R)\|^2_0\leq
C_{R,T}\cdot\sqrt{\eps}.\label{PP8}
\ee

On the other hand, setting
\ce
\cB(s,\eps,\eps'):=\sum_k\|B^\eps_k(s,u^\eps(s))-B^{\eps'}_k(s,u^{\eps'}(s))\|^2_0,
\de
by BDG's inequality and Young's inequality, we have
\ce
&&\mE\Bigg(\sup_{s\in[0,T\wedge\tau^{\eps,\eps'}_R]}|J_3(s)|\Bigg)\\
&\leq&C\mE\Bigg(\int^{T\wedge\tau^{\eps,\eps'}_R}_0\|v(s)\|^2_0\cdot
\cB(s,\eps,\eps')\dif s\Bigg)^{1/2}\\
&\leq&\frac{1}{2}\mE\Bigg(\sup_{s\in[0,T\wedge\tau^{\eps,\eps'}_R]}\|v(s)\|^2_0\Bigg)
+C\mE\Bigg(\int^{T\wedge\tau^{\eps,\eps'}_R}_0\cB(s,\eps,\eps')\dif s\Bigg).
\de
Thus, by (\ref{Io5})-(\ref{PP8}) and Lemma \ref{Le1} we obtain
\ce
\mE\Bigg(\sup_{s\in[0,T\wedge\tau^{\eps,\eps'}_R]}\|v(s)\|^2_0\Bigg)\leq C_{R,T}\cdot\sqrt{\eps}.
\de
Therefore, by Theorem \ref{Th1} and (\ref{Io6}) we get
\ce
&&\mE\left(\sup_{s\in[0,T]}\|v(s)\|^{2}_0\right)\\
&=&\mE\left(\sup_{s\in[0,T]}\|v(s)\|^{2}_0
\cdot 1_{\{\tau^{\eps,\eps'}_R\geq T\}}\right)
+\mE\left(\sup_{s\in[0,T]}\|v(s)\|^{2}_0\cdot 1_{\{\tau^{\eps,\eps'}_R<T\}}\right)\\
&\leq&\mE\Bigg(\sup_{s\in[0,T\wedge\tau^{\eps,\eps'}_R]}\|v(s)\|^{2}_0\Bigg)
+\left[\mE\left(\sup_{s\in[0,T]}\|v(s)\|^{4}_0\right)\right]^{1/2}
\cdot \left[\mP(\tau^{\eps,\eps'}_R<T)\right]^{1/2}\\
&\leq&C_{R,T}\cdot\sqrt{\eps}+C_T/R.
\de
Lastly, letting $\eps\downarrow 0$ and $R\rightarrow\infty$, yield the desired limit.
\end{proof}

{\it Proof of Theorem \ref{main}:}

First of all,
by Lemma \ref{Le7}, there is a $u(\cdot)\in L^2(\Omega, \cF,\mP; C([0,T];\mH^0))$ such that
\be
\lim_{\eps\rightarrow 0}\mE\left(\sup_{s\in[0,T]}\|u^\eps(s)-u(s)\|^{2}_0\right)=0,\label{Es6}
\ee
which together with Theorem \ref{Th1} yields that for any $p\geq 1$
$$
\mE\left(\sup_{t\in[0,T]}\|u(t)\|^{2p}_\cN\right)
+\int^T_0\mE\|u(s)\|^2_{\cN+1}\dif s\leq C_{p,T}.
$$

We now show that $u(t)$ satisfies (ii) of Theorem \ref{main}.
It suffices to prove that for any $v\in\mH^\infty$
\ce
\<v,u(t)\>_0&=&\<v,u_0\>_0+\int^t_0\<v,\cL u(s)\>_0\dif s+\int^t_0\<v,F(s,u(s))\>_0\dif s\\
&&+\sum_k\int^t_0\<v,B_k(s,u(s))\>_0\dif W^k(s), \ \ \mP-a.s..
\de
Note that
\ce
\<v,u^\eps(t)\>_0&=&\<v,u_0\>_0+\int^t_0\<v,\cT_\eps\cL \cT_\eps u^\eps(s)\>_0\dif s
+\int^t_0\<v,\cT_\eps F(s,\cT_\eps u^\eps(s))\>_0\dif s\\
&&+\sum_k\int^t_0\<v,\cT_\eps B_k(s,\cT_\eps u^\eps(s))\>_0\dif W^k(s), \ \ \mP-a.s..
\de
We only prove that the third term in the right hand side converges to the corresponding term,
that is, as $\eps\downarrow 0$
\be
\cP(t,\eps):=\int^t_0|\<v,\cT_\eps F(s,\cT_\eps u^\eps (s))-F(s,u(s))\>_0|\dif s
\stackrel{\text{\tiny $L^1(\Omega;\mP)$}}{\longrightarrow} 0.\label{Lim}
\ee

For any $R>0$, define the stopping time
$$
\tau_R^\eps:=\inf\{t>0: \|u^\eps(t)\|_N\wedge\|u(t)\|_N\geq R\}.
$$
Thus,
\be
\cP(t,\eps)=\cP(t,\eps)\cdot 1_{\{\tau^{\eps}_R\geq t\}}
+\cP(t,\eps)\cdot  1_{\{\tau^{\eps}_R<t\}}.\label{Po1}
\ee
For the first term of (\ref{Po1}), we have by {(\bf $\mathbf{H}1_N$})
\ce
&&\mE\big(\cP(t,\eps)\cdot  1_{\{\tau^{\eps}_R\geq t\}}\big)\leq
\mE(\cP(t\wedge\tau^{\eps}_R,\eps))\leq\\
&\leq&\mE\left(\int^{t\wedge\tau^{\eps}_R}_0|\<\cT_\eps v,
F(s,\cT_\eps u^\eps (s))-F(s,u(s))\>_0|\dif s\right)\\
&&+\mE\left(\int^{t\wedge\tau^{\eps}_R}_0|\<\cT_\eps v-v,
F(s,u(s))\>_0|\dif s\right)\\
&\leq&\|\cT_\eps v\|_1\cdot\mE\left(\int^{t\wedge\tau^{\eps}_R}_0
\|F(s,\cT_\eps u^\eps (s))-F(s,u(s))\|_{-1}\dif s\right)\\
&&+\|\cT_\eps v-v\|_1\cdot\mE\left(\int^{t\wedge\tau^{\eps}_R}_0
\|F(s,u(s))\|_{-1}\dif s\right)\\
&\leq&C_R\cdot \|v\|_1\cdot\int^{t}_0
\mE\|\cT_\eps u^\eps (s)-u(s)\|_0\dif s+C_{R,T}\cdot\eps\cdot\|v\|_3.
\de
For the second term of (\ref{Po1}), as above it is easy to see  by (\ref{Co7}) and
Theorem \ref{Th1} that
\ce
\mE(\cP(t,\eps)\cdot  1_{\{\tau^{\eps}_R<t\}})\leq C_T\cdot \mP(\tau^{\eps}_R<t)^{1/2}\leq C_T/R.
\de
Therefore, first letting $\eps\downarrow 0$ and then $R\rightarrow\infty$ for (\ref{Po1}) gives that
\ce
\lim_{\eps\downarrow 0}\mE\cP(t,\eps)=0.
\de
The uniqueness follows from similar calculations as in  proving Lemma \ref{Le7}.
The proof is thus complete.
\br
By (\ref{Es6}) and Theorem \ref{Th1}, using the interpolation inequality
(\ref{Int}) and H\"older's inequality, we in fact have for any $0<\a<\cN$
\ce
\lim_{\eps\downarrow 0}\mE\left(\sup_{s\in[0,T]}\|u^\eps(s)-u(s)\|^{2}_\a\right)=0.
\de
\er

\section{Strong Solutions of Semilinear SPDEs in Euclidean space}

Consider the following Cauchy problem of SPDE in $\mR^d$:
\ce
\left\{
\begin{array}{lcl}
\dif u(t,x)&=&\Big[\Delta u(t,x)+\sum_{i=1}^d \p_{x_i}f_i(t,\om,x,u(t,x))
+g(t,\om,x,u(t,x))\Big]\dif t\\
&&+\sum_k\sigma_k(t,\om,x,u(t,x))\dif W^k(t),\\
u(0,x)&=&u_0(x),
\end{array}
\right.
\de
where $\Delta:=\sum_{i=1}^d\p^2_{x_i}$ is the Laplace operator, and
the other coefficients are respectively measurable with respect to their variables:
\ce
&&f:\mR_+\times\Omega\times\mR^d\times\mR\to\mR^d,\\
&&g:\mR_+\times\Omega\times\mR^d\to\mR,\\
&&\sigma:\mR_+\times\Omega\times\mR^d\to l^2.
\de

We impose the following conditions on $f,g$ and $\sigma$:
\begin{enumerate}[{(\bf $\mathbf{A}1$)}]

\item For each $x\in \mR^d$, $z\in\mR$ and $t\geq 0$,
$f(t,\cdot,x,z), g(t,\cdot,x,z)$ and $\sigma(t,\cdot,x,z)$ are $\cF_t$-measurable.

\item  There exist $\kappa_1,\kappa_2>0$ and $h_0,h_1\in L^2(\mR)$ such that
for all $(t,\omega)\in\mR_+\times\Omega$, $x\in \mR^d$
and $z\in\mR$
$$
\|\p_zf(t,\om,x,z)\|_{\mR^d}+ |\p_z g(t,\om,x,z)|+\|\p_z \sigma(t,\omega,x,z)\|_{l^2}\leq \kappa_1,
$$
and for $j=0,1$
$$
\|\nabla_x^jf(t,\om,x,z)\|_{\mR^d}+ |\nabla_x^jg(t,\om,x,z)|
+\|\nabla_x^j\sigma(t,\omega,x,z)\|_{l^2}\leq \kappa_2|z|+h_j(x),
$$
where $\nabla_x=(\p_{x_1},\cdots,\p_{x_d})$ is the gradient operator.
\end{enumerate}

For $m\in\mN_0:=\{0\}\cup\mN$, let $\mW^{m}_2(\mR^d)$ be the usual Sobolev space in $\mR^d$, i.e.,
the completion of the space $C_0^\infty(\mR^d)$ of smooth functions with compact supports
with respect to the norm
\ce
\|u\|_m&:=&\left(\int_{\mR^d}|u(x)|^2\dif x+
\int_{\mR^d}|\nabla^mu(x)|^2\dif x\right)^{1/2}\\
&=&\left(\int_{\mR^d}|(I-\Delta)^{\frac{m}{2}}u(x)|^2\dif x\right)^{1/2}.
\de
We will set $\mH^m:=\mW^m_2(\mR^d)$ and $\cL:=\Delta$, and define
\be
F(t,\om,u)&:=&\sum_{i=1}^d \p_{x_i}f_i(t,\om,\cdot,u(\cdot))+g(t,\om,\cdot,u(\cdot)),\label{Oi1}\\
B_k(t,\om,u)&:=&\sigma_k(t,\om,\cdot,u(\cdot)),\ \ k\in\mN.\label{Oi2}
\ee
In the following, for the simplicity of notations, the variables $t$ and $\om$ in $F$
and $B$ will be hidden.
\bl\label{Le9}
Assume {(\bf $\mathbf{A}1$)} and {(\bf $\mathbf{A}2$)}. Then
$F$ and $B$ defined by (\ref{Oi1}) and (\ref{Oi2})
satisfy {(\bf $\mathbf{H}1_1$)} and {(\bf $\mathbf{H}2_1$)}.
\el
\begin{proof}
Notice that
$(I-\Delta)^{-\frac{1}{2}}\p_{x_i}$ and $(I-\Delta)^{-\frac{1}{2}}$ are
bounded linear operators from $L^2(\mR^d)$ to $L^2(\mR^d)$.
It is clear by {(\bf $\mathbf{A}2$)}  that for any $u,v\in\mH^1$
\ce
\|F(u)-F(v)\|^2_{-1}+\sum_k\|B_k(u)-B_k(v)\|_0^2&\leq& C\|u-v\|_0^2,\\
\|F(u)\|^2_0+\sum_k\|B_k(u)\|^2_0&\leq& C(\|u\|_1^2+1),
\de
and by integration by parts formula and Young's inequality
\ce
\<u,F(u)\>_0&\leq& \frac{1}{2}\|u\|^2_1+C\|f(\cdot,u(\cdot))\|_0^2
+C\|g(\cdot,u(\cdot))\|_0^2\\
&\leq& \frac{1}{2}\|u\|^2_1+C(\|u\|_0^2+1).
\de
Hence, {(\bf $\mathbf{H}1_1$)} holds.

For {(\bf $\mathbf{H}2_1$)}, we have
\ce
\<u,F(u)\>_1&\leq& \frac{1}{4}\|u\|^2_2+C\|F(u)\|_0^2\leq
\frac{1}{4}\|u\|^2_2+C(\|u\|_1^2+1)\leq\\
&\leq& \frac{1}{4}\|u\|^2_2+C(\|u\|_2\cdot\|u\|_0+1)\leq
\frac{1}{2}\|u\|^2_2+C(\|u\|_0^2+1),
\de
and by {(\bf $\mathbf{A}2$)}
\ce
\sum_k\|B_k(u)\|^2_1\leq C(\|u\|_1^2+1)\leq \delta\|u\|^2_2+C_\delta(\|u\|_0^2+1).
\de
The proof is complete.
\end{proof}

By Lemma \ref{Le9} and Theorem \ref{main}, we obtain the following result:
\bt\label{Th2}
Assume {(\bf $\mathbf{A}1$)} and {(\bf $\mathbf{A}2$)}.
For any $u_0\in\mH^1$, there exists a unique $\mH^1$-valued continuous and $\cF_t$-
adapted process $u(t)$ such that for any $T>0$ and $p\geq 1$
$$
\mE\left(\sup_{s\in[0,T]}\|u(s)\|^p_1\right)
+\int^T_0\mE\|u(s)\|^2_{2}\dif s<+\infty,
$$
and the following equation holds in $\mH^0$: for all $t\geq 0$
\ce
u(t,\cdot)&=&u_0(\cdot)+\int^t_0\Big[\Delta u(s,\cdot)+\sum_{i=1}^d \p_{x_i}f_i(s,\cdot,u(s,\cdot))
+g(s,\cdot,u(s,\cdot))\Big]\dif s\\
&&+\sum_k\int^t_0\sigma_k(s,\cdot,u(s,\cdot))\dif W^k(s),\ \ \mP-a.s..
\de
\et
\br
This result is not new (cf. \cite{MiRo}), however,
our general result can be used to treat the initial-boundary values problem
as follows.
\er
We now turn to the initial-boundary values problem. Let $\cO$ be a bounded smooth domain in  $\mR^d$.
Consider the following SPDE with Dirichlet boundary conditions:
\ce
\left\{
\begin{array}{lcl}
\dif u(t,x)&=&\Big[\Delta u(t,x)+\sum_{i=1}^d \p_{x_i}f_i(t,\om,x,u(t,x))
+g(t,\om,x,u(t,x))\Big]\dif t\\
&&+\sum_k\sigma_k(t,\om,x,u(t,x))\dif W^k(t),\\
u(t,x)&=&0,~(t,x)\in\mR_+\times\p\cO,\\
u(0,x)&=&u_0(x).
\end{array}
\right.
\de

For $m\in\mN_0$, let $\mW^m_2(\cO)$ and $\mW^{m,0}_2(\cO)$
be the usual Soblev spaces on $\cO$, which are
the respective completions of smooth functions spaces $C^\infty(\cO)$
and $C^\infty_0(\cO)$(with compact supports in $\cO$)
with respect to the norm:
\ce
\|f\|_{m,2}:=\left(\sum_{j=0}^m\int_\cO|\nabla^j f(x)|^2\dif x\right)^{1/2}.
\de

Set $\cL:=\Delta$ and $\sD(\cL):=\mW^2_2(\cO)\cap\mW^{1,0}_2(\cO)$.
Then $(\cL,\sD)$ forms a sectorial operator in $L^2(\cO)$(cf. \cite{Pa}),
and we have the corresponding $\mH^m$.
We remark that $\mH^1=\mW^{1,0}_2(\cO)$ and
$(I-\cL)^{-1/2}\p_{x_i}$ and $\p_{x_i}(I-\cL)^{-1/2}$ are  bounded linear operators in $L^2(\cO)$.

We assume that:
\begin{enumerate}[{(\bf $\mathbf{A}1'$)}]

\item For each $x\in \cO$, $z\in\mR$ and $t\geq 0$,
$f(t,\cdot,x,z), g(t,\cdot,x,z)$ and $\sigma(t,\cdot,x,z)$ are $\cF_t$-measurable.

\item  There exist $\kappa_1,\kappa_2>0$ and $h_0,h_1\in L^2(\cO)$ such that
for all $(t,\omega)\in\mR_+\times\Omega$, $x\in \cO$
and $z\in\mR$
$$
\|\p_zf(t,\om,x,z)\|_{\mR^d}+ |\p_z g(t,\om,x,z)|+\|\p_z \sigma(t,\omega,x,z)\|_{l^2}\leq \kappa_1,
$$
and for $j=0,1$
$$
\|\nabla_x^jf(t,\om,x,z)\|_{\mR^d}+ |\nabla_x^jg(t,\om,x,z)|
+\|\nabla_x^j\sigma(t,\omega,x,z)\|_{l^2}\leq \kappa_2|z|+h_j(x),
$$
\item One of the following conditions holds:
$$
\sigma(t,\om,x,0)=0\ \ \mbox{ or }\ \ \sigma(t,\om,x,z)=0\mbox{ for any $x\in\p\cO$}.
$$
\end{enumerate}
\br
Notice the following characterization of $\mW^{1,0}_2(\cO)$(cf. \cite{Ze}):
$$
\mbox{$u\in\mW^{1,0}_2(\cO)$
if and only if $u\in\mW^1_2(\cO)$ and $u(x)=0$ for almost all $x\in\p\cO$.}
$$
Thus,
{(\bf $\mathbf{A}2'$)} and {(\bf $\mathbf{A}3'$)} imply that if $u\in\mH^1$, then
$\sigma_k(t,\om,\cdot,u(\cdot))\in\mH^1$.
\er

Basing on the similar calculations as above, we have
\bt
Assume {(\bf $\mathbf{A}1'$)}-{(\bf $\mathbf{A}3'$)}.
For any $u_0\in\mH^1$, the same conclusions  as in Theorem \ref{Th2} hold.
Moreover, $u(t,x)=0$ for almost all $x\in\p\cO$ and any $t\geq 0$.
\et
\section{Stochastic Burgers and Ginzburg-Landau's equations on the real line}

In this section, we consider the following generalized stochastic Burgers
and Ginzburg-Landau equation on the real line $\mR$:
\ce
\dif u(t,x)&=&[\p^2_x u(t,x)+\p_x f(t,\om,u(t,x))+g(t,\om,x,u(t,x))]\dif t\\
&&+\sum_k \sigma_k(t,\om,x,u(t,x))\dif W^k(t),
\de
where the coefficients $f,g$ and $\sigma_k, k\in\mN$ are measurable with respect
to their variables, and satisfy the following assumptions:
\begin{enumerate}[{(\bf $\mathbf{B}$1)}]
\item For each $t\geq 0$ and $x,z\in\mR$, $f(t,\cdot,z),g(t,\cdot,x,z)$ and
$\sigma_k(t,\cdot,x,z), k\in\mN$ are $\cF_t$-measurable.
\item
For every $(t,\om)\in\mR_+\times\Omega$, $f(t,\om,\cdot)\in C^\infty(\mR)$,
and for each $m\in\mN$, there exist a $q_m\geq 0$ and $\kappa^f_{m}>0$
such that for all $(t,\om,z)\in\mR_+\times\Omega\times\mR$
$$
|\p^m_z f(t,\omega,z)|\leq \kappa^f_{m}(|z|^{q_m}+1),
$$
where $q_1<2$.
\item
For every $(t,\om)\in\mR_+\times\Omega$, $g(t,\om,\cdot,\cdot)\in C^\infty(\mR^2)$,
and for each $n\in\mN_0$ and $m\in\mN$,
there exist $l_{nm}, l_n\geq 1$, $\kappa^g_{nm},\kappa_n^g>0$
and $h^g_{n}\in L^2(\mR)$ such that for all $(t,\om,x,z)\in\mR_+\times\Omega\times\mR\times\mR$
\ce
|\p^n_x\p^m_z g(t,\omega,x,z)|&\leq& \kappa^g_{nm}\cdot(|z|^{l_{nm}}+1),\\
|\p^n_x g(t,\omega,x,z)|&\leq&\kappa^g_n\cdot|z|^{l_n}+h^g_n(x)
\de
where  $1\leq l_1<7$, and for some $\kappa^g>0$
$$
\p_zg(t,\omega,x,z)\leq\kappa^g.
$$

\item For every $(t,\om)\in\mR_+\times\Omega$ and $k\in\mN$,
$\sigma_k(t,\om,\cdot,\cdot)\in C^\infty(\mR^2)$,
and for each $n\in\mN_0$ and $m\in\mN$,
there exist $p_{nm}\geq 0, p_n\geq 1$, $\kappa^\sigma_{nm},\kappa_n^\sigma>0$
and $h^\sigma_{n}\in L^1(\mR)$ such that for all $(t,\om,x,z)\in\mR_+\times\Omega\times\mR\times\mR$
\ce
\sum_k|\p^n_x\p^m_z \sigma_k(t,\omega,x,z)|^2&\leq& \kappa^\sigma_{nm}\cdot(|z|^{p_{nm}}+1),\\
\sum_k|\p^n_x \sigma_k(t,\omega,x,z)|^2&\leq&\kappa^\sigma_n\cdot|z|^{2p_n}+h^\sigma_n(x),
\de
where $p_{01}=0$ and $1\leq p_1<5$.
\end{enumerate}
\br\label{R1}
The condition $\p_zg(t,\omega,x,z)\leq \kappa^g$ implies that
$$
z\cdot g(t,\omega,x,z)\leq\kappa^g\cdot z^2-g(t,\omega,x,0)\cdot|z|,\ \ \forall z\in\mR.
$$
In particular, $f(z)=z^2$ and $g(z)=z-z^{2n-1}$
satisfy {(\bf $\mathbf{B}$2)} and {(\bf $\mathbf{B}$3)}.
\er

Let $\mW^m_2(\mR)$ be the usual Sobolev spaces on $\mR$. We need the following Gagliardo-Nirenberg
inequality: for any $p\in[2,+\infty]$, $m\in\mN$ and $u\in \mW^m_2(\mR)$
 (cf. \cite[p.24, Theorem 9.3]{Fr})
\be
\|u\|_{L^p}\leq C\|u\|_m^{\frac{p-2}{2mp}}\|u\|_0^{\frac{2mp-p+2}{2mp}}\leq C\|u\|_m.\label{Ep}
\ee
In the following, we take $\mH^m=\mW^m_2(\mR)$ and $\cL=\p_x^2$, and define for $u\in L^2(\mR)$
\be
F(t,\omega,u)&:=&\p_x f(t,\omega,u(\cdot))+g(t,\omega,\cdot,u(\cdot)),\label{DF1}\\
B_k(t,\omega,u)&:=&\sigma_k(t,\omega,\cdot,u(\cdot)),\ \ k\in\mN.\label{DF2}
\ee
As in the previous section, we hide the variables $t$ and $\om$ without confusions.
We have
\bl
Assume {(\bf $\mathbf{B}$1)}-{(\bf $\mathbf{B}$4)}, then $F$ and $B$
defined by (\ref{DF1}) and (\ref{DF2}) satisfy {(\bf $\mathbf{H}1_1$)}.
\el
\begin{proof}
Noting that
$(I-\p_x^2)^{-\frac{1}{2}}\p_x$ and $(I-\p_x^2)^{-\frac{1}{2}}$ are
bounded linear operators from $L^2(\mR)$ to $L^2(\mR)$
and the following elementary formula
\be
\phi(u)-\phi(v)=(u-v)\int^1_0\phi'(s(u-v)+v)\dif s,\label{For}
\ee
by {(\bf $\mathbf{B}$2)} and {(\bf $\mathbf{B}$3)} we have for any $u,v\in\mH^1$
\ce
\|F(u)-F(v)\|_{-1}
&\leq& C\|f(u(\cdot))-f(v(\cdot))\|_0+\|g(\cdot,u(\cdot))-g(\cdot,v(\cdot))\|_0\\
&\leq& C(\|u\|^{q_1\vee l_{01}}_{L^\infty}+\|v\|^{q_1\vee l_{01}}_{L^\infty}+1)\cdot\|u-v\|_0\\
&\leq& C(\|u\|^{q_1\vee l_{01}}_1+\|v\|^{q_1\vee l_{01}}_1+1)\cdot\|u-v\|_0.
\de
It is obvious by {(\bf $\mathbf{B}$4)} with $p_{01}=0$ and (\ref{For}) that
\ce
\sum_k\|B_k(u)-B_k(v)\|_0^2\leq 2\kappa^\sigma_{01}\cdot\|u-v\|^2_0.
\de

Moreover, by integration by parts formula we have for $u\in\mH^1$
\ce
&&\int_\mR u(x)\p_x f(u(\cdot))\dif x=-\int_\mR \p_x u(x)f(u(x))\dif x\\
&=&-\int_\mR \p_x\left(\int^{u(\cdot)}_0f(r)\dif r\right)\dif x=0.
\de
By Remark \ref{R1}, we have
\ce
\int_\mR u(x)g(x,u(x))\dif x\leq C(\|u\|^2_0+1).
\de
Hence
$$
\<u,F(u)\>_0\leq C(\|u\|^2_0+1).
$$

On the other hand, as above it is easy to see that for some $p>1$
\ce
\|F(u)\|_0&\leq&C\|u\|_{L^\infty}^{q_1}(\|u\|_1+1)+C\|u\|^{l_0-1}_{L^\infty}(\|u\|_0+1)\\
&\leq& C(\|u\|^{p}_1+1)
\de
and
\ce
\sum_k\|B_k(u)\|^2_0\leq C(\|u\|^2_0+1).
\de
The proof is complete.
\end{proof}
For verifying {(\bf $\mathbf{H}2_\cN$}), we need the following elementary differential formula,
which can be proved by induction method.
\bl\label{Le8}
Let $\phi\in C^\infty(\mR^2)$. For $m\geq 3$ and $u\in\mH^m$, we have
\ce
\p^{m}_x \phi(\cdot,u)&=&(\p_u\phi)(\cdot,u)\p^{m}_xu+
m[(\p^2_u\phi)(\cdot,u)\p_x u+(\p_x\p_u\phi)(\cdot,u)]\cdot\p^{m-1}_xu\\
&&+a(\p_x^{m-2} u,\cdots,\p_x u)+(\p^m_x\phi)(\cdot,u),
\de
where $a$ is a multi-function and has at most polynomial growth with respect to its variables.
\el
\bl\label{Le10}
For any $m\in\mN$ and $u\in\mH^m$, we have
\ce
\<\p^m_xu,\p^{m}_xg(\cdot,u)\>_0+\|\p^m_x f(u)\|_0^2\leq
\frac{1}{4}\|u\|_{m+1}^2+C(\|u\|_{m-1}^{\a_m}+1).
\de
\el
\begin{proof}
For $m=1$, we have by {(\bf $\mathbf{B}$3)} and (\ref{Ep})
\ce
\<\p_xu,\p_xg(\cdot,u)\>_0&=&\<\p_xu,(\p_xg)(x,u)\>_0+
\<\p_xu,(\p_ug)(x,u)\p_x u\>_0\\
&\leq& \|\p_xu\|_0\cdot(\kappa^g_1\cdot\||u|^{l_1}\|_0+\|h^g_1\|_0)
+\kappa^g\cdot\|\p_x u\|^2_0\\
&\leq& C(\|u\|^2_1+1)+\kappa^g_1\cdot\|u\|_1\cdot\|u\|^{l_1}_{L^{2l_1}}\\
&\leq& C\cdot(\|u\|_2\cdot\|u\|_0+1)+C\|u\|^{\frac{l_1+1}{4}}_2
\cdot\|u\|_0^{\frac{3l_1+5}{4}}
\de
and by {(\bf $\mathbf{B}$2)}
\ce
\|\p_x f(u)\|_0^2&\leq& (\kappa^f_1)^2\cdot\int_\mR(|u|^{q_1}+1)^2\cdot|\p_x u|^2\dif x\\
&\leq& (\kappa^f_1)^2\cdot(\|u\|_{L^\infty}^{2q_1}+1)\cdot\|u\|^2_1\\
&\leq& C\|u\|_2^{\frac{q_1}{2}+1}\|u\|_0^{\frac{3q_1}{2}+1}+C\cdot\|u\|_2\cdot\|u\|_0.
\de
Since $l_1<7$ and $q_1<2$, by Young's inequality we get for some $\a_1>1$
\ce
\<\p_xu,\p_xg(u)\>_0+\|\p_x f(u)\|_0^2\leq \frac{1}{4}\|u\|_2^2+C(\|u\|_0^{\a_1}+1).
\de

For $m=2$,  noticing that
\ce
\p^2_xg(\cdot,u)&=&(\p_ug)(x,u)\p^2_xu
+(\p^2_ug)(x,u)(\p_xu)^2\\
&&+2(\p_x\p_ug)(x,u)\p_xu+(\p^2_xg)(x,u),
\de
we have by {(\bf $\mathbf{B}$3)}  and (\ref{Ep})
\ce
\<\p^2_xu,\p^2_xg(\cdot,u)\>_0&\leq&\kappa^g\cdot\|\p^2_x u\|^2_0+
\|\p^2_xu\|_0\cdot\Big[\|(\p^2_ug)(\cdot,u)(\p_xu)^2\|_0\\
&&+2\|(\p_x\p_ug)(\cdot,u)\p_xu\|_0
+\|(\p^2_xg)(\cdot,u)\|_0\Big]\\
&\leq&\kappa^g\cdot\|u\|^2_2+
\|u\|_2\cdot\Big[\kappa^g_{02}(\|u\|^{l_{02}}_{L^\infty}+1)\|\p_xu\|_{L^4}^4\\
&&+2\kappa^g_{11}(\|u\|^{l_{11}}_{L^\infty}+1)\|\p_xu\|_0
+\kappa^g_2\|u\|^{l_2}_{L^{2l_2}}+\|h^g_2\|_0\Big]\\
&\leq&\kappa^g\cdot\|u\|^2_2+
\|u\|_2\cdot\Big[C(\|u\|^{l_{02}}_1+1)\|u\|_{2}\cdot\|u\|_1^3\\
&&+C(\|u\|^{l_{11}}_1+1)\|u\|_1
+\kappa^g_2\|u\|^{l_2}_1+\|h^g_2\|_0\Big],\\
&\leq&C\|u\|_3\|u\|_1\Big[1+C(\|u\|^{l_{02}}_1+1)\cdot\|u\|_1^3\Big]\\
&&+C(\|u\|^{2(l_{11}+1)}_1+\|u\|^{2l_{2}}_1+1),
\de
and by {(\bf $\mathbf{B}$2)}
\ce
\|\p^2_x f(u)\|_0^2&\leq& 2\int_\mR|(\p_uf)(u)\p^2_x u|^2\dif x+
2\int_\mR|(\p^2_uf)(u)(\p_x u)^2|^2\dif x\\
&\leq& C(\|u\|^{q_1}_{L^\infty}+1)\cdot\|\p^2_x u\|^2_0+
C(\|u\|^{q_2}_{L^\infty}+1)\cdot\|\p_x u\|^4_{L^4}\\
&\leq& C(\|u\|^{q_1+1}_{1}+1)\cdot\|u\|_3+
C(\|u\|^{q_2+7/2}_1+1)\cdot\|u\|^{1/2}_3.
\de
Hence, for some $\a_2>1$
\ce
\<\p^2_xu,\p^2_xg(\cdot,u)\>_0+\|\p^2_x f(u)\|_0^2
\leq \frac{1}{4}\|u\|_3^2+C(\|u\|_1^{\a_2}+1).
\de
The higher derivatives can be estimated similarly by Lemma \ref{Le8}.
\end{proof}
\bl
Assume {(\bf $\mathbf{B}$1)}-{(\bf $\mathbf{B}$4)}, then $F$ and $B$
defined by (\ref{DF1}) and (\ref{DF2}) satisfy {(\bf $\mathbf{H}2_\cN$})
for any $\cN\in\mN$.
\el
\begin{proof}
For any $m\in\mN$, we have by Lemma \ref{Le10}
\ce
\<u,F(u)\>_m&=&\<u,F(u)\>_0+\<\p^m_xu,\p^m_xF(u)\>_0\\
&\leq& C(\|u\|^2_0+1)+\<\p^m_xu,\p^{m+1}_xf(u)\>_0
+\<\p^m_xu,\p^{m}_xg(\cdot,u)\>_0\\
&\leq& C(\|u\|^2_0+1)+\frac{1}{4}\|u\|^2_{m+1}+\|\p^{m}_xf(u)\|^2_0
+\<\p^m_xu,\p^{m}_xg(\cdot,u)\>_0\\
&\leq& \frac{1}{2}\|u\|_{m+1}^2+C(\|u\|_{m-1}^{\a_m}+1).
\de
It is similar to the proof of Lemma \ref{Le10}
that for any $\delta>0$ and some $\beta_m\geq 1$
$$
\sum_k\|B_k(u)\|_m^2\leq \delta\|u\|^2_{m+1}+C(\|u\|^{\beta_m}_{m-1}+1).
$$
The proof is complete.
\end{proof}
Summarizing the above calculations, we obtain by Theorem \ref{main}
\bt
Under {(\bf $\mathbf{B}$1)}-{(\bf $\mathbf{B}$4)},
for any $u_0\in\mH^\infty$, there exists a unique process
$u(t)\in\mH^\infty$ such that
\begin{enumerate}[(i)]
\item For any $m\in\mN$, the process $t\mapsto u(t)\in\mH^m$ is $\cF_t$-adapted and continuous, and
for any $T>0$ and $p\geq 2$
$$
\mE\left(\sup_{s\in[0,T]}\|u(s)\|^p_m\right)<+\infty.
$$
\item For almost all $\om$ and all $t\geq 0$, $x\in\mR$
\ce
u(t,x)&=&u_0(x)+\int^t_0[\p^2_xu(s,x)+\p_x f(s,u(s,x))+g(s,x,u(s,x))]\dif s\\
&&+\sum_k\int^t_0\sigma_k(s,x,u(s,x))\dif W^k(s).
\de
\end{enumerate}
\et

\section{Stochastic tamed 3D Navier-Stokes equation in $\mR^3$}

In this and next sections,
we shall use bold-face letters $\u=(u_1,u_2,u_3),\cdots,$
to denote the velocity fields in $\mR^3$(or $\mR^2$).

Consider the following stochastic tamed 3D Navier-Stokes equation with viscosity constant
$\nu=1$ in $\mR^3$:
\be
\dif \u(t)&=&\Big[\Delta \u(t)-(\u(t)\cdot\nabla)\u(t)+\nabla p(t)-g_N(|\u(t)|^2)\u(t)\Big]\dif t\no\\
&&+\sum_{k=1}^\infty\Big[\nabla \tilde p_k(t)
+\h_k(t,\om,x,\u(t))\Big]\dif W^k_t\label{Ns1}
\ee
subject to the incompressibility condition
\be
\div(\u(t))=0,\label{In}
\ee
and the initial condition
\be
\u(0)=\u_0,\label{Ns2}
\ee
where $p(t,x)$ and $\tilde p_k(t,x)$ are unknown scalar functions,
 $N>0$ and the taming function $g_N:\mR^+\mapsto\mR^+$ is a smooth function such that
\be
\label{Con}\left\{
\begin{array}{ll}
g_N(r):=0, \quad r\in[0,N], \\
g_N(r):=r-N,\quad  r\geq N+2,\\
0\leq g'_N(r)\leq C,\quad  r\geq 0,\\
|g^{(k)}_N(r)|\leq C_{k},\quad  r\geq 0, k\in\mN,
\end{array}
\right.
\ee
and $\h_k, k\in\mN$ satisfy that
\begin{enumerate}[{(\bf $\mathbf{C1}$)}]
\item For each $k\in\mN$ and $t\geq0 $, $x,z\in \mR^3$,
$\h_k(t,\cdot,x,z)$ are $\cF_t$-measurable.
\item For every $(t,\om)\in\mR_+\times\Omega$ and $k\in\mN$,
$\h_k(t,\om,\cdot,\cdot)\in C^\infty(\mR^3\times\mR^3;\mR^3)$,
and for each $n\in\mN_0$ and $m\in\mN$,
there exist  $\kappa_{nm},\kappa_n>0$
and $b_{n}\in L^1(\mR^3)$ such that for all
$(t,\om,x,z)\in\mR_+\times\Omega\times\mR^3\times\mR^3$
$$
\sum_k|\p^n_x\p^m_z \h_k(t,\omega,x,z)|^2\leq \kappa_{nm},
$$
and
$$
\sum_k|\p^n_x \h_k(t,\omega,x,z)|^2\leq\kappa_n\cdot|z|^2+b_n(x).
$$
\end{enumerate}

For $m\in\mN_0$, set
\be
\mH^m:=\{\u\in \mW^{m}_2(\mR^3)^3: \div(\u)=0\},\label{Div}
\ee
where $\div$ is taken in the sense of Schwartz distributions.
The following Gagliardo-Nirenberg
inequality will be used frequently below  (cf. \cite[p.24, Theorem 9.3]{Fr}):
for $r\in[2,+\infty]$ and $3(\frac{1}{2}-\frac{1}{r})\leq m\in\mN$
\be
\|\u\|_{L^r}^r\leq C\|\u\|^{\frac{3(r-2)}{2m}}_{m}\|\u\|_0^{r-\frac{3(r-2)}{2m}}.\label{Sob}
\ee

Let $\sP$ be the orthogonal projection operator from $L^2(\mR^3)^3$ to $\mH^0$.
It is well known
that $\sP$ can be restricted to a bounded linear operator from $\mW^{m}_2(\mR^3)^3$ to $\mH^m$,
and that $\sP$ commutes with the
derivative operators (cf. \cite{Lions}). For any $\u\in\mH^0$ and $\v\in L^2(\mR^3)^3$, we have
$$
\<\u,\v\>_{\mH^0}:=\<\u,\sP\v\>_{\mH^0}=\<\u,\v\>_{L^2}.
$$

For $\u\in\mH^2$, define
\be
\cL\u&:=&\sP\Delta \u,\label{FF0}\\
F(\u)&:=&-\sP((\u\cdot\nabla)\u)-\sP(g_N(|\u|^2)\u),\label{FF1}\\
B_k(t,\om,\u)&:=&\sP(\h_k(t,\om,\cdot,\u)).\label{FF2}
\ee

Using $\sP$ to act on both sides of Eq.(\ref{Ns1}), we may consider the following
equivalent equation of (\ref{Ns1})-(\ref{Ns2}):
\ce
\dif \u(t)=\Big[\cL\u(t)+F(\u(t))\Big]\dif t
+\sum_{k=1}^\infty B_k(t,\u(t))\dif W^k_t, \ \ \u(0)=\u_0.
\de
\bl\label{Le01}
Under {(\bf $\mathbf{C2}$)}, the operators $F$ and $B$ defined by (\ref{FF1})
and (\ref{FF2}) satisfy   {(\bf $\mathbf{H}1_2$)}.
\el
\begin{proof}
Noting that
$(I-\cL)^{-\frac{1}{2}}\nabla \sP$ and $(I-\cL)^{-\frac{1}{2}}\sP$ are
bounded linear operators from $\mH^0$ to $\mH^0$ and
$$
\sP((\u\cdot\nabla)\u)=\sum_{i=1}^3\p_{x_i} \sP(u_i\cdot \u),
$$
we have by (\ref{Sob})
\ce
&&\|F(\u)-F(\v)\|_{-1}\\
&\leq& C\|\u\otimes\u-\v\otimes\v\|_0
+C\|g_N(|\u|^2)\u-g_N(|\v|^2)\v\|_0\\
&\leq&C(\|\u\|_{L^\infty}+\|\v\|_{L^\infty})\|\u-\v\|_0
+C(\|\u\|^2_{L^\infty}+\|\v\|^2_{L^\infty})\|\u-\v\|_0\\
&\leq&C(\|\u\|^2_{2}+\|\v\|^2_{2}+1)\|\u-\v\|_0.
\de
Moreover, it is obvious by {(\bf $\mathbf{C2}$)}  that
$$
\sum_k\|B_k(t,\u)-B_k(t,\u)\|^2_0\leq C\|\u-\v\|_0^2
$$
and
$$
\sum_k\|B_k(t,\u)\|^2_0\leq C(\|\u\|_0^2+1).
$$

On the other hand, observing that
$$
\<\u,(\u\cdot\nabla)\u\>_0=\frac{1}{2}\<\u,\nabla|\u|^2\>_0=0
$$
and
$$
\<\u,g_N(|\u|^2)\u\>_0\geq 0,
$$
we have
\be
\<\u,F(\u)\>_0\leq 0.\label{F0}
\ee

Lastly, it is easy to see by (\ref{Sob}) that
\ce
\|F(\u)\|_0&\leq& \|\u\|_{L^\infty}\|\u\|_1+\|\u\|_{L^6}^3\\
&\leq&C\|\u\|^{3/4}_2\cdot\|\u\|_0^{1/4}\cdot\|\u\|_1+C\|\u\|_1^3\\
&\leq&C(\|\u\|_2+\|\u\|_1^5+1).
\de
The proof is complete.
\end{proof}

In order to check {(\bf $\mathbf{H}2_\cN$}), we need
the following Moser-type calculus inequality
 (cf. \cite[p. 294, Proposition 21.77]{Ze}):
\bl\label{Le02}
For any $m\in\mN_0$, there exists a $C_m>0$ such that for any $\u,\v\in L^\infty\cap\mH^m$
\be
\|\u\cdot\v\|_{m}\leq C_m\Big[\|\u\|_{L^\infty}\cdot\|\v\|_m+\|\v\|_{L^\infty}\cdot\|\u\|_m\Big].\label{Cal}
\ee
\el

We now prove the following key estimate.
\bl\label{Le3}
Under {(\bf $\mathbf{C2}$)}, the operators $F$ and $B$ defined by (\ref{FF1})
and (\ref{FF2}) satisfy   {(\bf $\mathbf{H}2_\cN$}) for any $\cN\in\mN$.
\el
\begin{proof}
Let us first verify (\ref{Co4}) for $m=1$.
By Young's inequality, we have
\ce
-\<\u,\sP((\u\cdot\nabla)\u)\>_1&\leq&\frac{1}{4}\|\u\|^2_2
+\|(\u\cdot\nabla)\u\|^2_0\\
&\leq&\frac{1}{4}\|\u\|^2_2+\||\u|\cdot|\nabla\u|\|^2_0,
\de
where
\ce
|\u|^2=\sum_{k=1}^3|u_k|^2,\quad |\nabla\u|^2=\sum_{k,i=1}^3|\p_{x_i}u_k|^2.
\de

From the expression of $g_N$, we  also have
\ce
&&-\<\u,\sP(g_N(|\u|^2)\u)\>_1\\
&=&-\<\nabla (g_N(|\u|^2)\u), \nabla\u\>_0
-\<g_N(|\u|^2)\u, \u\>_0\\
&=&-\sum_{k,i=1}^3\int_{\mR^3} \p_{x_i}u_k\cdot \p_{x_i} (g_N(|\u|^2)u_k)\dif x
-\int_{\mR^3}|\u|^2\cdot g_N(|\u|^2)\dif x\\
&\leq&-\sum_{k,i=1}^3\int_{\mR^3} \p_{x_i} u_k\cdot \left(g_N(|\u|^2)\cdot \p_{x_i}u_k
-g'_N(|\u|^2)\p_{x_i} |\u|^2\cdot u_k\right)\dif x\\\
&=&-\int_{\mR^3}|\nabla \u|^2\cdot g_N(|\u|^2)\dif x
-\frac{1}{2}\int_{\mR^3}g'_N(|\u|^2)|\nabla|\u|^2|^2\dif x\\
&\leq&-\int_{\mR^3}|\nabla \u|^2\cdot |\u|^2\dif x+N\|\nabla\u\|^2_0.
\de
Hence
\be
\<\u,F(\u)\>_1&\leq& \frac{1}{4}\|\u\|^2_2+N\|\u\|^2_1\label{F1}\\
&\leq&\frac{1}{2}\|\u\|^2_2+C_N\|\u\|^2_0.\no
\ee

For $m\geq 2$, by the calculus inequality (\ref{Cal})
\ce
-\<\u,\sP((\u\cdot\nabla)\u)\>_m&=&-\<(I-\Delta)^{m/2}\u,(I-\Delta)^{m/2}((\u\cdot\nabla)\u)\>_0\\
&\leq&\frac{1}{8}\|\u\|_{m+1}^2+2\|(\u\cdot\nabla)\u\|^2_{m-1}\\
&\leq&\frac{1}{8}\|\u\|_{m+1}^2+
C_m(\|\u\|^2_{L^\infty}\|\u\|^2_m+\|\nabla\u\|^2_{L^\infty}\|\u\|^2_{m-1}).
\de
Noting that by  Agmon's inequality (cf. \cite{Kre}),
$$
\|\u\|^2_{L^\infty}\leq C\|\u\|_2\cdot\|\u\|_1,
$$
we have
\ce
\|\u\|^2_{L^\infty}\|\u\|^2_m\leq C\|\u\|^3_m
\cdot\|\u\|_1\leq C_{m}\|\u\|^{\frac{3}{2}}_{m+1}
\cdot\|\u\|^{\frac{3}{2}}_{m-1}\cdot\|\u\|_1
\de
and
\ce
\|\nabla\u\|^2_{L^\infty}\|\u\|^2_{m-1}&\leq& C\|\u\|_3\cdot\|\u\|_2
\cdot\|\u\|^2_{m-1}\\
&\leq& C\|\u\|^{\frac{3}{2}}_3\cdot\|\u\|_1^{\frac{1}{2}}
\cdot\|\u\|^2_{m-1}\\
&\leq& C\|\u\|^{\frac{3}{2}}_{m+1}\cdot\|\u\|_{m-1}^{\frac{5}{2}}.
\de
Thus, we get by Young's inequality
$$
-\<\u,\sP((\u\cdot\nabla)\u)\>_m\leq \frac{1}{4}\|\u\|_{m+1}^2+C_m\|\u\|^{10}_{m-1}.
$$

Let us now look at the term $-\<\u,\sP(g_N(|\u|^2)\u)\>_m$.
By the calculus inequality (\ref{Cal}) again, we have
\ce
\|g_N(|\u|^2)\u\|^2_{m-1}
\leq C_m(\|g_N(|\u|^2)\|^2_{L^\infty}\cdot\|\u\|^2_{m-1}
+\|g_N(|\u|^2)\|^2_{m-1}\cdot\|\u\|^2_{L^\infty})
\de
and
\ce
\||\u|^2\|^2_{m-1}\leq C_m\|\u\|^2_{m-1}\cdot\|\u\|^2_{L^\infty}.
\de
Noting that for any $k\geq 2$
$$
|g^{(k)}_N(r)|=0\ \ \mbox{ on $r<N$ and $r>N+2$},
$$
we have
$$
\|g_N(|\u|^2)\|_{m-1}\leq C_m\||\u|^2\|_{m-1}+C_{N,m}(\|\u\|^{\gamma_m}_{m-2}+1),\ \ \gamma_m>2.
$$
As above, we obtain
\ce
&&-\<\u,\sP(g_N(|\u|^2)\u)\>_m\\
&\leq&\frac{1}{8}\|\u\|_{m+1}^2+2\|g_N(|\u|^2)\u\|^2_{m-1}\\
&\leq&\frac{1}{8}\|\u\|_{m+1}^2+C_m\|\u\|^4_{L^\infty}\cdot\|\u\|^2_{m-1}
+C_{N,m}\|\u\|^2_{L^\infty}\cdot(\|\u\|^{\gamma_m}_{m-2}+1)\\
&\leq&\frac{1}{4}\|\u\|_{m+1}^2+C_m(\|\u\|^{\gamma_m}_{m-1}+1),\ \ \gamma_m>2.
\de
Combining the above calculations yields that for some $\a_m\geq 1$
\ce
\<\u,F(\u)\>_m\leq\frac{1}{2}\|\u\|_{m+1}^2+C_m(\|\u\|^{\a_m}_{m-1}+1).
\de

We now check (\ref{Co6}). For $m=1$, we have by {(\bf $\mathbf{NS}$)} and (\ref{Int})
\be
\sum_{k=1}^\infty \|\nabla B_k(t,\u)\|_0^2&\leq&\sum_{k=1}^\infty \|(\nabla_x\h_k)(t,\u)\|_0^2
+\sum_{k=1}^\infty \|(\p_u\h_k)(t,\u)\nabla\u\|_0^2\no\\
&\leq&C(\|\u\|_0^2+1)+C\|\nabla\u\|_0^2\label{F4}\\
&\leq&C(\|\u\|_0^2+1)+C\|\u\|_2\cdot\|\u\|_0\no\\
&\leq&\delta\|\u\|^2_2+C(\|\u\|^2_0+1).\no
\ee
The higher derivatives can be calculated similarly.
The proof is complete.
\end{proof}

Thus, we obtain the following main result in this section.
\bt\label{Main}
Under {\bf ($\mathbf{C1}$)} and {\bf ($\mathbf{C2}$)}, for any $\u_0\in\mH^\infty$
there exists a unique solution $\u(t)\in\mH^\infty$ to Eq.(\ref{Ns1}) such that
\begin{enumerate}[(i)]
\item For any $m\in\mN$, the process $t\mapsto u(t)\in\mH^m$ is $\cF_t$-adapted and continuous, and
for any $T>0$ and $p\geq 2$
$$
\mE\left(\sup_{s\in[0,T]}\|\u(s)\|^p_m\right)<+\infty.
$$
\item For almost all $\om$ and all $t\geq 0$,
\ce
\u(t)&=&\u_0+\int^t_0[\Delta \u(s)-\sP((\u(s)\cdot\nabla)\u(s))]\dif s\\
&&-\int^t_0\sP(g_N(|\u(s)|^2)\u(s))\dif s+\sum_k\int^t_0\sP(\h_k(s,\u(s)))\dif W^k(s).
\de
\end{enumerate}
\et

\section{Stochastic 2D Navier-Stokes equation in $\mR^2$}

Consider the following stochastic 2D Navier-Stokes equation in $\mR^2$
\be
\dif \u(t)&=&\Big[\Delta \u(t)-(\u(t)\cdot\nabla)\u(t)+\nabla p(t)\Big]\dif t\no\\
&&+\sum_{k=1}^\infty\Big[\nabla \tilde p_k(t)
+\h_k(t,\om,x,\u(t))\Big]\dif W^k_t\label{Ns01}
\ee
subject to the incompressibility condition $\div(\u(t))=0$,
and the initial condition $\u(0)=\u_0$. As in Section 6,
the functions $p(t,x)$ and $\tilde p_k(t,x)$ are unknown scalar functions,
$\u(t,x)$ is the velocity field in $\mR^2$, and $\h_k, k\in\mN$ satisfy
{(\bf $\mathbf{C1}$)} and {(\bf $\mathbf{C2}$)}.

For $m\in\mN_0$, set
\ce
\mH^m:=\{\u\in \mW^{m}_2(\mR^2)^2: \div(\u)=0\}.
\de
We also have the projection operator $\sP$ from $L^2(\mR^2)^2$ to $\mH^0$.

Our main result in this section is that
\bt\label{bt1}
Assume {(\bf $\mathbf{C1}$)} and {(\bf $\mathbf{C2}$)}. For any $\u_0\in\mH^\infty$,
there exists a unique solution $\u(t)\in\mH^\infty$ to Eq.(\ref{Ns01}) such that
\begin{enumerate}[(i)]
\item
For any $m\in\mN$, the process $t\mapsto u(t)\in\mH^m$ is $\cF_t$-adapted and continuous,
and for some sequence of stopping times $\tau_n\uparrow\infty$
and any $T>0$ and $p\geq 2$
$$
\mE\left(\sup_{s\in[0,T\wedge\tau_n]}\|\u(s)\|^p_m\right)\leq C_{n,T}.
$$
\item For almost all $\om$ and all $t\geq 0$,
\ce
\u(t)=\u_0+\int^t_0[\Delta \u(s)-\sP((\u(s)\cdot\nabla)\u(s))]\dif s
+\sum_k\int^t_0\sP(\h_k(s,\u(s)))\dif W^k(s).
\de
\end{enumerate}
\et

We can not directly use Theorem \ref{main} to prove this result because $\mH^1$ does not embedded in
$L^\infty(\mR^2)^2$. We shall first consider the modified equation like (\ref{Ns1}), then use
stopping times technique to obtain the existence of smooth solutions for equation (\ref{Ns01}).

Note that the result in Section 6 also holds for 2D Navier-Stokes equation.
In what follows, we shall use the same notations as in Section 6.
The solution in Theorem \ref{Main} corresponding to the tamed function $g_N$
is denoted by $\u_N(t)$.

Before proving Theorem \ref{bt1}, we prepare the following lemma.
\bl\label{Lp1}
There exists a constant $C>0$ independent of $N$ such that for any $\u\in\mH^3$,
$$
2\<\u,F(\u)\>_2+\sum_k\|B_k(s,\u)\|^2_2
\leq \|\u\|^2_3+C\|\u\|^2_2\cdot(1+\|\u\|^2_1)\cdot (1+\|\u\|^2_0)+C.
$$
\el
\begin{proof}
By Young's inequality, we have
\ce
\<\u,F(\u)\>_2\leq \frac{1}{4}\|\u\|^2_3+2\|(\u\cdot\nabla)\u\|_1^2+2\|g_N(|\u|^2)\u\|_1^2.
\de
Noticing that
$$
\|\u\|^4_{L^4}\leq C\|\u\|^2_1\cdot\|\u\|^2_0,\quad
\|\u\|^4_{L^\infty}\leq C\|\u\|^2_2\cdot\|\u\|^2_0,
$$
we have
\ce
\|\nabla((\u\cdot\nabla)\u)\|_0^2&\leq& 2\||\u|\cdot|\nabla^2\u|\|^2_0+2\|\nabla\u\|^4_{L^4}\\
&\leq& 2\|\u\|^2_{L^4}\cdot\|\nabla^2\u\|^2_{L^4}+C\|\u\|^2_{2}\cdot\|\u\|^2_{1}\\
&\leq& C\|\u\|_{0}\cdot\|\u\|_{1}\cdot\|\u\|_{2}\cdot\|\u\|_3+C\|\u\|^2_{2}\cdot\|\u\|^2_{1}\\
&\leq& \frac{1}{4}\|\u\|^2_3+C\|\u\|^2_2\cdot\|\u\|^2_1\cdot (1+\|\u\|^2_0)
\de
and
$$
\|\nabla(g_N(|\u|^2)\u)\|_0^2\leq C\|\u\|^4_{L^\infty}\cdot\|\u\|^2_1
\leq C\|\u\|^2_2\cdot\|\u\|^2_0\cdot \|\u\|^2_1.
$$

Moreover, it is easy to see by {(\bf $\mathbf{C2}$)} that
\ce
\sum_k\|B_k(s,\u)\|^2_2\leq C\|\u\|^2_{2}\cdot(\|\u\|^2_{1}+1)+C.
\de
The result now follows.
\end{proof}

Now we can give the proof of Theorem \ref{bt1}.

\begin{proof}
By Ito's formula,   we have for any $p\geq 1$
\be
\|\u_N(t)\|^{2p}_m=\|\u_0\|^{2p}_m+\sum_{j=1}^5I^{(p)}_{mj}(t),\label{Io11}
\ee
where
\ce
I^{(p)}_{m1}(t)&:=&2p\int^t_0\|\u_N(s)\|^{2(p-1)}_m\<\u_N(s),\Delta\u_N(s)\>_m\dif s\\
I^{(p)}_{m2}(t)&:=&2p\int^t_0\|\u_N(s)\|^{2(p-1)}_m\<\u_N(s),F(\u_N(s))\>_m\dif s\\
I^{(p)}_{m3}(t)&:=&p\sum_k\int^t_0\|\u_N(s)\|^{2(p-1)}_m\|B_k(s,\u_N(s))\|_m^2\dif s\\
I^{(p)}_{m4}(t)&:=&2p\sum_{k=1}^\infty\int^t_0\|\u_N(s)\|^{2(p-1)}_m\<\u_N(s),
B_k(s,\u_N(s))\>_m\dif W^k_s\\
I^{(p)}_{m5}(t)&:=&2p(p-1)\sum_{k=1}^\infty\int^t_0\|\u_N(s)\|^{2(p-2)}_m|\<\u_N(s),
B_k(s,\u_N(s))\>_m|^2\dif s.
\de

For $m=0$,
by (\ref{F0}) and taking expectations for (\ref{Io11}), it is easy to see that
for any $T>0$ and $p\geq 1$
\be
\sup_{t\in[0,T]}\mE\|\u_N(t)\|^{2p}_0+\int^T_0\mE\Big(\|\u_N(s)\|^{2(p-1)}_0\cdot
\|\u_N(s)\|^{2}_1\Big)\dif s\leq C_{T,p}.\label{F2}
\ee
Here and after, the constant $C_{T,p}$ is independent of $N$.

Define
$$
\lambda_N(t):=\int^t_0(1+\|\u_N(s)\|^2_1)\cdot(1+\|\u_N(s)\|^2_0)\dif s
$$
and stopping times
$$
\theta_t^N:=\inf\{s\geq 0: \lambda_N(s)\geq t\}.
$$
Then $\lambda_N^{-1}(t)=\theta_t^N$ and $\theta^N_t\leq t$.
Moreover, by (\ref{F2})  we have for any $T>0$
\be
\lim_{t\rightarrow\infty}\sup_N \mP(\theta^N_t<T)=0.\label{Op1}
\ee

From (\ref{Io11}) and using Lemma \ref{Lp1}, we have for $m=2$ and $p=1$
\ce
\|\u_N(\theta^N_t)\|^2_2&\leq&\|\u_0\|^{2}_2-\int^{\theta^N_t}_0\|\u_N(s)\|^2_3\dif s
+C\int^{\theta^N_t}_0\|\u_N(s)\|^2_2~
\dif \lambda_N(s)\\
&&+C\theta^N_t+I^{(1)}_{24}(\theta^N_t)\\
&\leq&\|\u_0\|^{2}_2-\int^{\theta^N_t}_0\|\u_N(s)\|^2_3\dif s+C\int^t_0\|\u_N(\theta^N_s)\|^2_2~
\dif s\\
&&+Ct+I^{(1)}_{24}(\theta^N_t).
\de
Taking expectations and using Gronwall's inequality yield that for any $T>0$
\ce
\sup_{t\in[0,T]}\mE\|\u_N(\theta^N_t)\|^2_2+\mE\left(\int^{\theta^N_T}_0\|\u_N(s)\|^2_3\dif s
\right)\leq C_T(\|\u_0\|^{2}_2+1).
\de
Using the same trick as in proving (\ref{Es2}), we further have
\ce
\mE\left(\sup_{t\in[0,\theta^N_T]}\|\u_N(t)\|^2_2\right)
=\mE\left(\sup_{t\in[0,T]}\|\u_N(\theta^N_t)\|^2_2\right)\leq C_T(\|\u_0\|^{2}_2+1).
\de

Set for $N>1$
$$
\tau_N:=\inf\left\{t\geq 0: \sup_{x\in\mR^2}|\u_N(t,x)|\geq\sqrt{N}\right\}.
$$
Then for any $T>0$
\be
\lim_{N\rightarrow\infty}\mP(\tau_N\leq \theta_T^N)&=&
\lim_{N\rightarrow\infty}\mP\left(\sup_{t\in[0,\theta_T^N]}\sup_{x\in\mR^2}|\u_N(t,x)|\geq\sqrt{N}\right)\no\\
&\leq&\lim_{N\rightarrow\infty}\mP\left(\sup_{t\in[0,\theta_T^N]}\|\u_N(t)\|_2\geq\sqrt{N}\right)\no\\
&\leq&\lim_{N\rightarrow\infty}\mE\left(\sup_{t\in[0,\theta_T^N]}\|\u_N(t)\|_2^2\right)/N\no\\
&\leq&\lim_{N\rightarrow\infty}C_T(\|\u_0\|^{2}_2+1)/N=0.\label{Op2}
\ee

On the other hand, it is clear that $\u_N$ satisfy Eq.(\ref{Ns01}) on $[0,\tau_N]$.
By the uniqueness of solution to Eq.(\ref{Ns01}), we know that
$$
\tau_N\leq\tau_{N+1}, \ \ \mP-a.s..
$$
Noting that for any $M,N,T>0$
\ce
\mP(\tau_N<M)&=&\mP(\tau_N<M; \theta_T^N<M)+\mP(\tau_N<M; \theta_T^N\geq M)\\
&\leq&\mP(\theta_T^N<M)+\mP(\tau_N\leq \theta_T^N),
\de
we have by (\ref{Op1}) and (\ref{Op2})
$$
\lim_{N\rightarrow\infty}\mP(\tau_N<M)=0,
$$
which means that
$$
\lim_{N\rightarrow\infty}\tau_N=\infty, \ \ \mP-a.s..
$$
The whole proof is complete.
\end{proof}

{\bf Acknowledgements:}

Part of this work was done while the author was a fellow of Alexander-Humboldt Foundation.
He would  like to thank Professor Michael R\"ockner, his host in Bielefeld University,  for
providing him  a stimulating environment.

\end{document}